\documentclass{article}
\usepackage{graphicx} % Required for inserting images
\usepackage[english]{babel} % Langue
\usepackage{geometry}
\usepackage{csquotes}
\usepackage{appendix}
\usepackage{amsthm}
\usepackage{mathrsfs}
\usepackage{amsmath}
\usepackage{tikz-cd}
\usepackage{amssymb}
\usepackage{stmaryrd}
\usepackage[T1]{fontenc}
\usepackage[all]{xy}
\usepackage{enumitem}
\newlist{tirets}{enumerate}{1}
\setlist[tirets,1]{label=-}

\theoremstyle{definition}
\newtheorem{definition}{Definition}[section]
\theoremstyle{plain}
\newtheorem{theorem}[definition]{Theorem}
\newtheorem{proposition}[definition]{Proposition}
\newtheorem{lemma}[definition]{Lemma}
\newtheorem{corollary}[definition]{Corollary}
\theoremstyle{remark}
\newtheorem{remark}[definition]{Remark}

\title{A Homological Criterion for the Almost-Existence of Hamiltonian Chords}
\author{Antoine Rodrigues}
\date{October 2025}

\begin{document}

\maketitle

\begin{abstract}
    We establish a criterion on wrapped Floer homology of an exact Lagrangian submanifold in a Liouville domain, which ensures the almost-existence of Hamiltonian chords near a given energy level. To this purpose we introduce a relative version of the coisotropic Hofer-Zehnder capacity, and compare it with some other capacities that are defined using filtered wrapped Floer homology. Our main application is the almost-existence of Hamiltonian chords near any hypersurface if the wrapped Floer homology vanishes.
\end{abstract}

\tableofcontents

\section{Introduction}

\subsection{Motivation}

The existence of Hamiltonian or Reeb orbits is one of the central topics in symplectic and contact geometry. One of the questions in this area is Arnold’s chord conjecture. It asserts the existence of Reeb chords for any Legendrian submanifold of a contact manifold and any contact form, and it has only been proved in this generality in dimension 3, in the paper \cite{ConjCordesDim3} by Hutchings and Taubes.\\

This is related to the existence of Hamiltonian chords for a Lagrangian submanifold in symplectic geometry: If $(W,\theta)$ is a Liouville domain, its boundary $\partial W$ has a canonical contact form $\alpha = \theta|_{\partial W}$, and a Lagrangian submanifold $L$ that is invariant under the Liouville flow $Z$ near $\partial W$ gives a Legendrian submanifold $L\cap\partial W$ of $\partial W$. If $H$ is an autonomous Hamiltonian that only depends on the coordinate associated to the Liouville vector field near $\partial W$, then the chords of $H$ that are contained in the hypersurface $\partial W$ correspond to Reeb chords for the Legendrian submanifold $L\cap\partial W$. In this context, a recent work \cite{ConjCordesConormal} by Broćić, Cant and Shelukhin proves the chord conjecture in the specific case where the Lagrangian is a conormal bundle and the Liouville domain is the unit cotangent bundle of a compact manifold. Ritter also shows that this conjecture holds on all contact manifolds that are the boundary of a Liouville domain whose symplectic homology vanishes, in particular the boundaries of subcritical Stein manifolds (\cite{Ritter}, Example 11.2).\\

The goal in this paper is to obtain existence results for Hamiltonian chords with boundary on some Lagrangian submanifold of a Liouville domain, near a given energy level. The strategy we employ follows the same line as the one developed in \cite{dg-floer2} (Section 5) for the search for periodic Hamiltonian orbits. It is based on the observation that, while the question of the existence of chords on a hypersurface may be difficult to settle, it is possible to obtain \emph{almost-existence} results. More precisely, we define the almost-existence property in the neighborhood of a hypersurface:

\begin{definition}
\label{almost-existence}
    Given a Hamiltonian $H$ on a symplectic manifold $W$ and a Lagrangian submanifold $L$, we say that the almost-existence property for Hamiltonian chords (resp. periodic orbits) holds for a hypersurface $\Sigma = H^{-1}\{a\}$ if there exists a neighborhood $(a-\epsilon,a+\epsilon)$ of $a$ such that, for almost all values $b \in (a-\epsilon,a+\epsilon)$ in the sense of the Lebesgue measure, the hypersurface $H^{-1}\{b\}$ contains a Hamiltonian chord (resp. a periodic orbit).
\end{definition}

A classical result states that this almost-existence property, for periodic orbits instead of chords, can be guaranteed by the finiteness of the Hofer–Zehnder capacity (see \cite[Chapter 4, Theorem 4]{HoferZehnder} or \cite{Struwe}), and this can in turn be ensured by a criterion involving the symplectic homology of the Liouville domain $W$ (this idea appears in \cite{Irie}). The authors of \cite{dg-floer2} then obtain the following homological criterion for the almost-existence of periodic orbits:

\begin{theorem} \cite[Theorem 5.1]{dg-floer2}
    \label{critère homologique => presque existence}
    Let $\Sigma$ be a hypersurface bounding a relatively compact open subset $U$. Suppose there exists a homology class $\alpha \in H_*(W,\partial W)$ such that $i_W(\alpha)=0$ and $j_{U!}(\alpha)\neq 0$. Then the almost-existence property for contractible periodic orbits holds in a neighborhood of $\Sigma$.
\end{theorem}

Here, $i_W$ is the canonical morphism between Morse homology and symplectic homology, and $j_{U!} : H_*(W,\partial W)\to H_*(U,\partial U)$ is the shriek morphism in Morse homology induced by the inclusion $j_U : U\hookrightarrow W$. This theorem was established for a generalized version of symplectic homology with differential graded coefficients, which allows the assumptions of the theorem to be satisfied in many more situations than with classical symplectic homology. The authors of \cite{dg-floer2} thus generalize Morse homology and Floer homology to incorporate differential graded coefficients.\\

We will follow here the approach proposed in \cite{dg-floer2} to analyze the Lagrangian case, without considering differential graded coefficients. In the following, Morse and Floer homologies will be taken with coefficients in $\mathbb{Z}/2\mathbb{Z}$.

\subsection{Results}

Let $L$ be an exact Lagrangian submanifold of a Liouville domain $(W,\theta)$ such that the intersection $\partial L = \partial W\cap L$ is transverse, and $\theta|_L$ vanishes in a neighborhood of $\partial L$. The following theorem, analogous to Theorem \ref{critère homologique => presque existence} in the Lagrangian case, is the main result of this paper:

\begin{theorem}
\label{critère homologique => cordes}
    Let $\Sigma$ be a hypersurface bounding a relatively compact open subset $U$, intersecting $L$ transversely. Suppose there exists a homology class $\alpha\in H_*(L,\partial L)$ such that $i_L(\alpha)=0$ and $j_{U!}(\alpha)\neq0$. Then the almost-existence property for contractible Hamiltonian chords holds in a neighborhood of $\Sigma$.
\end{theorem}

Here, $i_L$ is the canonical morphism between Morse homology $H_*(L,\partial L)$ and wrapped Floer homology $SH_{*,0}(W,L)$ (see part \ref{Morse -> symplectique}), and $j_{U!}$ is the shriek map induced by the inclusion $L\cap U \hookrightarrow L$.

\begin{remark}
\label{Reduction to capacities ?}
    As we will see in part \ref{Le critère homologique}, the proof of Theorem \ref{critère homologique => cordes} goes in the following way:
    $$\text{Homological criterion for }\Sigma=\partial U \Longrightarrow c_{HZ}^0(W,\overline{U};L)<\infty \Longrightarrow \text{Almost-existence near }\Sigma$$
    where $c_{HZ}^0(W,\overline{U};L)$ is the $\pi_1$-sensitive relative Lagrangian capacity that we define in part \ref{Capacités et orbites}. If the homological criterion is valid for every open set $U$ (for example when $SH_*(W,L)=0$), then the absolute capacity $c_{HZ}^0(W;L)$ is finite. These considerations are the same in the case of periodic orbits, so we have the following implications:
    \[\begin{tikzcd}
        & SH_*(W)=0 \arrow[d,Rightarrow]\arrow[r,Rightarrow] & c_{HZ}^0(W)<\infty \\
        & SH_*(W,L)=0 \arrow[r,Rightarrow] & c_{HZ}^0(W;L)<\infty.
    \end{tikzcd}\]
    Therefore, a natural question to ask, to which we do not give an answer in this paper, is whether the finiteness of the capacity $c_{HZ}^0(W)$ implies the finiteness of the Lagrangian capacity $c_{HZ}^0(W;L)$.\\
\end{remark}

An immediate consequence of Theorem \ref{critère homologique => cordes} is the following:

\begin{corollary}
\label{SH(W,L)=0=>chords}
    If $SH_*(W,L)=0$, then the almost existence property for contractible Hamiltonian chords is satisfied for any hypersurface that is transverse to $L$ and of the form $\Sigma=\partial U$, with $U$ a relatively compact open set.
\end{corollary}

\textit{Proof.} If $SH_*(W,L)=0$, then the homological criterion of Theorem \ref{critère homologique => cordes} is satisfied for all $U$ with the class $\alpha=[L,\partial L]$. Indeed, we have $j_{U!}[L,\partial L] = [U\cap L,\partial U\cap L]\neq0$.

\qed

\begin{remark}
    As we will see, the homological criterion can be stated in a form that involves only the wrapped Floer homology that corresponds to contractible chords $SH_{*,0}(W,L)$. However, the vanishing of this last homology is equivalent to the vanishing of $SH_*(W,L)$ since $SH_*(W,L)$ has a ring structure, with unit $i_L[L,\partial L]\in SH_{*,0}(W,L)$ (see \cite[Theorem 6.14]{Ritter}).\\
\end{remark}

We will now give some applications of Corollary \ref{SH(W,L)=0=>chords}.\\

\textbf{The case $\mathbf{SH_*(W)=0}$.} According to Ritter \cite[Theorem 6.17]{Ritter}, $SH_*(W,L)$ carries a natural $SH_*(W)$-module structure, where $SH_*(W)$ denotes the symplectic homology of $W$. It then suffices that this latter homology vanishes in order to conclude that $SH_*(W,L)=0$.\\

This is the case, for instance, of the symplectic ball $\overline{B}^n(0,1)\subset\mathbb{C}^n$ (see \cite{ReliqueFloer1} or \cite{Oancea_homologie_symplectique}, Section~3), or more generally of any compact convex domain in $\mathbb{C}^n$. A more general class of Liouville domains with vanishing symplectic homology is that of \emph{flexible Weinstein domains} (see \cite[Definition 11.29]{Stein2Weinstein}). A proof of this fact can be found in \cite[Theorem 3.2]{MurphySiegel18} (see also \cite{BourgeoisEkholmEliashberg12}). Furthermore, a version of Künneth formula for symplectic homology \cite{Oan06} shows that, if $W_1,W_2$ are Liouville domains such that $SH_*(W_1)=0$, then $SH_*(W_1\times W_2)=0$. Therefore, we get the following corollary of Theorem \ref{critère homologique => cordes}:

\begin{corollary}
\label{corollaire SH=0}
    Assume that $W$ is a compact convex domain in $\mathbb{C}^n$, a flexible Weinstein domain or a product of such with another Liouville domain. Then the almost existence property for contractible Hamiltonian chords is satisfied for any hypersurface that is transverse to $L$ and of the form $\Sigma=\partial U$, with $U$ a relatively compact open set.
\end{corollary}

Observe that this result applies in particular to subcritical Weinstein domains, since they are always flexible. Alternatively, it suffices to notice that any subcritical Weinstein manifold is symplectomorphic to $(\mathbb{C}\times V,\omega_0\oplus\omega_1)$, with $V$ a Stein manifold and $\omega_0$ the standard symplectic form on $\mathbb{C}$ by \cite[Theorem 14.16]{Stein2Weinstein}, so Corollary \ref{corollaire SH=0} also applies.\\

\textbf{Lagrangians that are displaceable from the boundary.} A subset $A$ of a symplectic manifold $W$ is called \emph{displaceable} from another subset $B$ if there exists a compactly supported Hamiltonian isotopy $H : [0,1]\times W\to\mathbb{R}$ such that $\phi^1_H(A)\cap B=\emptyset$. In \cite{KimKimKwon22}, the authors prove that, if the Liouville domain $W$ (or, equivalently, its boundary $\partial W$) is displaceable from $\hat{L}$ in the symplectic completion $\hat{W}$, then the wrapped Floer homology vanishes \cite[Corollary 2.15]{KimKimKwon22}. Actually, this result gives a more quantitative information about the Hofer-Zehnder capacity we use (see Remark \ref{Reduction to capacities ?}). Indeed, Theorem 2.14 in \cite{KimKimKwon22} shows that, for any Morse homology class $\alpha\in H_*(L,\partial L)$, we have
\begin{equation}
\label{c(a)<e}
    c(\alpha)\leq e(W;L)
\end{equation}
where $c(\alpha)$ is the quantity we define in the beginning of part \ref{Le critère homologique}, and
$$e(W;L) := \inf\{||K||, K\in C^\infty_c([0,1]\times W), \phi_K^1(\partial W)\cap\hat{L}=\emptyset\}$$
is the \emph{displacement energy} (here $||.||$ stands for the Hofer norm). Considering the capacity
$$c_{SH}(W;L) = \sup\{c(\alpha),\alpha\in H_*(L,\partial L)\} = c([L,\partial L])$$
as in Definition \ref{SH capacity}, and combining the inequality (\ref{c(a)<e}) with Theorem \ref{c_HZ < c_s < c_SH}, we get
$$c_{HZ}^0(W;L)\leq c_\sigma(W;L)\leq c_{SH}(W;L) \leq  e(W;L).$$

The displaceability condition can be ensured if the skeleton
$$\text{Sk}(W):=\bigcap_{t>0}\phi_Z^{-t}(W)$$
is disjoint from $L$. Indeed, in that case, since both $L$ and $W$ are compact, there exists $T>0$ such that $L\cap\phi_Z^{-T}(W)=\emptyset$, and one can construct a compactly-supported \emph{Hamiltonian} flow that acts like $\phi_Z$ on $L$.\footnote{
Let us describe the local form of this Hamiltonian flow. Let $L_t=\phi_Z^{-t}(L)$. For $t$ small enough, $L_t$ is contained in a Weinstein neighborhood of $L$, identified with a neighborhood of the zero section in $T^*L$, and $L_t$ is the graph of some 1-form $\beta_t$ on $L$. Then one can show that the exactness of $L$ implies that $\beta_t=df_t$ for some functions $f_t$ on $L$. Therefore, the flow of the Hamiltonian isotopy $H^t = f_t\circ\pi$
sends $L$ on $L_t$. The desired Hamiltonian flow can then be obtained by composition of a finite number of those flows.}
This proves the following corollary:

\begin{corollary}
    Assume that $W$ is a Weinstein domain such that the skeleton of $W$ is disjoint from $L$. Then the almost existence property for contractible Hamiltonian chords is satisfied for any hypersurface that is transverse to $L$ and of the form $\Sigma=\partial U$, with $U$ a relatively compact open set.
\end{corollary}

Once again, we recover the vanishing of wrapped Floer homology if $W$ is a flexible Weinstein domain: in that case, the skeleton has dimension at most $n-1$, so possibly after displacing it by a small Hamiltonian isotopy, it will not intersect $L$.\\

\begin{remark}
    \label{PE => E}
    We say that the hypersurface $\Sigma$ is \emph{of contact type} if there exists another vector field $Z'$ such that $\mathscr{L}_{Z'}\omega=\omega$, defined on an open neighborhood of $\Sigma$, that is transverse to $\Sigma$. In this case, the form $\theta':= i_{Z'}\omega$ is a primitive of $\omega$, and its restriction $\alpha'=\theta'|_\Sigma$ is a contact form.
    
    The almost-existence property from Definition \ref{almost-existence} implies the existence of chords (or periodic orbits) on $\Sigma$ if the hypersurface is of contact type. Indeed, for a contact type hypersurface, the flow of the vector field $Z'$ yields a diffeomorphism $V\simeq(-\epsilon,\epsilon)\times\Sigma$, where $V$ is an open neighborhood of $\Sigma$. Moreover, denoting by $\rho$ the coordinate on $(-\epsilon,\epsilon)$ and $r=e^\rho$, the symplectic structure induced on $(-\epsilon,\epsilon)\times\Sigma$ is given by $d(r\alpha')$ and we have $Z'=\partial_\rho$. Setting
\begin{align*}
    H : (-\epsilon,\epsilon)\times\Sigma & \longrightarrow \mathbb{R}\\
    H(\rho,x) & = e^\rho=r
\end{align*}
we observe that $X_H(\rho,x) = rR(x)$, where $R$ is the Reeb vector field on $\Sigma$. In particular, $[X_H,Z']=0$, so $Z'$ preserves the Hamiltonian trajectories contained in $V$. If $\Sigma$ satisfies the almost-existence property for Hamiltonian chords (resp. periodic orbits), then there exists a chord (resp. a periodic orbit) on a hypersurface $\Sigma_\rho = \{\rho\}\times\Sigma$ with $\rho\in(-\epsilon,\epsilon)$, which can then be transported via $Z'$ to obtain a chord (resp. a periodic orbit) on $\Sigma$.
\end{remark}

In particular, for a contact type hypersurface, the homological criterion in Theorem \ref{critère homologique => cordes} implies directly the existence of chords on $\Sigma$. We can now deduce an application of Theorem \ref{critère homologique => cordes} to contact geometry. As explained in the motivation, Arnold's chord conjecture asserts the existence of Reeb chords on any contact manifold, for any Legendrian submanifold and any contact form.

\begin{definition}
    Let $M$ be a contact manifold, and $N$ be a Legendrian submanifold of $M$. We say that the pair $(M,N)$ is \emph{exact fillable} if $M$ is the boundary of some Liouville domain $(W,\theta)$, its contact structure is given by the contact form $\theta|_{\partial W}$, and there is an exact Lagrangian submanifold $L$ of $W$ such that $\theta|_L=0$ near $\partial L$, and $N=\partial L$.
\end{definition}

Let $(M,N)$ be such a pair. Since $M=\partial W$ is a contact type hypersurface, if the homological criterion in Theorem \ref{critère homologique => cordes} is satisfied for $(W,L)$ and $U=W$ (that is, $i_L$ is not injective), we obtain the existence of Reeb chords on $(M,N)$ for the contact form $\alpha=\theta|_L$. But the Arnold's chords conjecture states the existence of Reeb chords for any contact form that represents the same contact structure. Let $\beta$ be another contact form that defines the same contact structure on $\partial W$. Then we have $\beta = f\alpha$ for some positive function $f$ on $\partial W$, and the graph of $f$ defines a hypersurface $\Sigma_f$ of the symplectic completion $\hat{W}$. This hypersurface is of contact type, equipped with the contact form $\theta|_{\Sigma_f}=f\alpha=\beta$, so if $H$ is a Hamiltonian such that $\Sigma_f$ is an energy level of $H$, then the Hamiltonian chords on $\Sigma_f$ are in bijection with the Reeb chords for this new contact form $\beta$. Therefore, we have the following corollary:

\begin{corollary}
    Let $M$ be a contact manifold and $N$ be a Legendrian submanifold of $M$ such that $(M,N)$ is exact fillable. Assume that the filling $(W,L)$ is such that
    $$i_L : H_*(L,\partial L) \longrightarrow SH_{*,0}(W,L)$$
    is not injective. Then Arnold chord conjecture holds for $(M,N)$: for any contact form that represents the contact structure of $M$, there is at least one Reeb chord in $M$ with endpoints on $N$.
\end{corollary}

\subsection{Structure of the paper}

As it was already the case for periodic orbits, the almost existence property for Hamiltonian chords will come from the finiteness of a capacity of Hofer-Zehnder flavor. Such a capacity was introduced by Lisi and Rieser in \cite{c_coisotrope}, and part \ref{Capacités et orbites} is dedicated to defining a relative version of it, inspired by the work of Ginzburg and Gürel \cite{c_relative}. We also show that the finiteness of this capacity implies a certain almost-existence property for Hamiltonian chords. In part \ref{Convention and overview}, we recall the main constructions for the Floer theory that we use, and fix some conventions. Part \ref{Le critère homologique} is dedicated to the proof of Theorem \ref{critère homologique => cordes}. The strategy we adopt is to bound the Hofer-Zehnder capacity we defined before, using operations on Floer homology.

\section{A relative version of the Lagrangian Hofer-Zehnder capacity}
\label{Capacités et orbites}

In \cite{c_coisotrope}, Lisi and Rieser introduce a capacity $c(W;N,\sim)$ for any coisotropic submanifold $N$ and any equivalence relation $\sim$ that is coarser than the one given by the isotropic foliation. They obtain an almost-existence result for Hamiltonian chords with endpoints that are on the same leaf whenever this capacity is finite \cite[Theorem 4.2]{c_coisotrope}. The case of a Lagrangian submanifold is a particular instance of this.\\

In \cite{c_relative}, Ginzburg and Gürel define a relative capacity $c(W,Z)$ associated to a symplectic manifold $W$ and a compact subset $Z$ \cite[Definition 2.8]{c_relative}. The finiteness of this capacity implies an almost existence property for periodic orbits near a given hypersurface that bounds a domain containing $Z$ \cite[Theorem 2.14]{c_relative}.\\

Our goal is to combine these two approaches to define a relative Lagrangian Hofer-Zehnder capacity, and derive the same kind of almost-existence result. Let $(W,\omega)$ be a symplectic manifold, $Z$ a compact subset of $W$, and $L$ a Lagrangian submanifold of $W$.

\begin{definition}
\label{H(W,Z;L)}
A Hamiltonian $H : W\to\mathbb{R}$ is called \emph{admissible} (for the Hofer-Zehnder capacity) if it satisfies the following two conditions:
\begin{enumerate}
\item $H=0=\max(H)$ on $W\backslash K$ where $K$ is a compact subset such that both $K$ and $W\backslash K$ intersect $L$.
\item $H=\min(H)$ on an open neighborhood $U$ of $Z$ such that both $U$ and $W\backslash U$ intersect $L$.
\end{enumerate}
We denote by $\mathcal{H}(W,Z;L)$ the set of such Hamiltonians.
\end{definition}

\begin{definition}
Let $x : \mathbb{R}_+\to W$ satisfy $x(0)\in L$. We define the return time of $x$ by:
$$\tau(x) = \inf\{t>0, x(t)\in L\}\in [0,+\infty].$$
\end{definition}

\begin{definition}
\label{c_HZ(W,Z;L)}
The \emph{relative Lagrangian Hofer-Zehnder capacity} is defined as
$$c_{HZ}(W,Z;L) = \sup\{-\min(H), H\in\underline{\mathcal{H}}(W,Z;L)\}$$
where $\underline{\mathcal{H}}(W,Z;L)$ is the set of Hamiltonians $H\in\mathcal{H}(W,Z;L)$ such that $H$ admits no non-constant chord $x$ with $x(0)\in L$ and $\tau(x)\leq1$.
\end{definition}

\begin{theorem}
\label{c(W,Z;L) => cords}
Let $H$ be a proper Hamiltonian such that $H=\min(H)$ on $Z$, and let $\lambda$ be a regular value of $H$ such that the hypersurface $H^{-1}\{\lambda\}$ intersects the submanifold $L$ transversely. Assume that $c_{HZ}(W,Z;L)<\infty$. Then for every $\epsilon>0$, there exists a regular value $\mu\in[\lambda-\epsilon,\lambda+\epsilon]$ such that $H^{-1}\{\mu\}$ contains a chord $x$ between two points of $L$ with $0<\tau(x)<\infty$.
\end{theorem}

\textit{Proof.} We combine the ideas of the proofs of Theorem 2.4 in \cite{c_relative} and Theorem 4.2 in \cite{c_coisotrope} to obtain the desired result. Given $H,\lambda,\epsilon$ satisfying the assumptions of the theorem, set $m=\min(H)$. Since $m$ is never a regular value of $H$, we have $\lambda>m$, and by possibly reducing $\epsilon$, we may assume that $\lambda-m\geq\epsilon$. Since the transversality condition is open, we may also assume that for every $\mu\in[\lambda-2\epsilon,\lambda+2\epsilon]$, $\mu$ is a regular value of $H$ and $L$ intersects the hypersurface $H^{-1}\{\mu\}$ transversely.\\

We construct an admissible Hamiltonian $F$ as follows: let $f : \mathbb{R}\to\mathbb{R}$ be a smooth function such that
\begin{align*}
    f(s) & = -c_{HZ}(W,Z;L)-1 & s\leq \lambda-\epsilon, \\
    f(s) & = 0 & s \geq \lambda,\\
    f'(s) & >0 & \lambda-\epsilon<s<\lambda.
\end{align*}
We then set $F = f\circ H$, and one can check that the Hamiltonian $F$ is admissible in the sense of Definition \ref{H(W,Z;L)}, when considering $K = H^{-1}[m,\lambda]$ and $U = H^{-1}(-\infty,\lambda-\epsilon)$.\\

Since $-\min(F)>c_{HZ}(W,Z;L)$, there exists a non-constant trajectory $x$ of $F$ such that $x(0)\in L$ and $\tau(x)\leq 1$. This orbit satisfies
$$\frac{dx}{dt} = X_F(x(t)) = f'(H(x(t)))X_H(x(t)),$$
and therefore it is contained in an energy level of $H$:
$$\frac{d}{dt}H(x(t)) = dH\left(\frac{dx}{dt}\right) = -\omega(X_H(x),X_F(x)) = 0.$$
Let $\mu = H(x(t))$ and $T = f'(\mu)$. We then have $\tfrac{dx}{dt} = T X_H(x)$, with $T\neq0$ since $x$ is non-constant, and by setting $\tilde{x}(t) = x(t/T)$ we obtain a non-constant chord that is contained in the energy level $H^{-1}\{\mu\}$. The fact that $f'(\mu)\neq0$ implies that $\mu\in(\lambda-\epsilon,\lambda)\subset[\lambda-\epsilon,\lambda+\epsilon]$. To prove that $\tilde{x}$ is the desired chord it remains to rule out the case $\tau(\tilde{x})=0$, and this follows from the transversality of $L$ and $H^{-1}\{\mu\}$: let $y\in L$ satisfying the transversality condition
$$T_yL + \text{Ker}(d_yH) = T_yW,$$
and such that $X_H(y)\in T_yL$. Then,
\begin{align*}
    & \forall v\in T_yL,\quad \omega(X_H(y),v)\in\omega(T_yL,T_yL)=0,\\
    & \forall v\in\text{Ker}(d_yH),\quad \omega(X_H(y),v) = -d_yH(v)=0,
\end{align*}
so $X_H(y)=0$. Therefore, since the trajectory $\tilde{x}$ is non-constant, it can never be tangent to the Lagrangian $L$, which shows that $\tau(\tilde{x})>0$.

\qed

Furthermore, this density result can be improved to obtain an almost existence property, as it was the case for the classical Hofer-Zehnder capacity (see \cite{Struwe}, or \cite[Chapter 4, Theorem 4]{HoferZehnder}).

\begin{theorem}
\label{c(W,Z;L) => presque existence}
Let $H$ be a proper Hamiltonian such that $H=\min(H)$ on $Z$, and let $\lambda$ be a regular value of $H$ such that the hypersurface $H^{-1}\{\lambda\}$ intersects the submanifold $L$ transversely. Assume that $c_{HZ}(W,Z;L)<\infty$. Then there exists $\epsilon>0$ such that for almost every $\mu\in[\lambda-\epsilon,\lambda+\epsilon]$, the hypersurface $H^{-1}\{\mu\}$ contains a chord $x$ with $0<\tau(x)<\infty$.
\end{theorem}

The proof of this fact is based on the following idea: in our situation, each energy level $H^{-1}\{\mu\}$ bounds the relatively compact open set $B_\mu = H^{-1}(\min(H)-\epsilon,\mu)$, and if we assume that the function
$$c : \mu \mapsto c_{HZ}(B_\mu,Z;L)$$
is Lipschitz on some interval $[\lambda,\lambda+\epsilon)$, then one can show that there exists a Hamiltonian chord that is contained in $H^{-1}\{\lambda\}$ (see \cite[Chapter 4, Theorem 3]{HoferZehnder} for a proof in the case of periodic orbits). Moreover, since the function $c$ is monotone, it is almost-everywhere differentiable by a theorem of Lebesgue, and therefore the Lipschitz condition above is satisfied for almost every $\lambda$, which completes the argument.

\begin{remark}
\label{pi_1 sensitive capacities}
    Let now $\eta$ be a homotopy class in $\pi_1(W,L)$. We define a capacity $c_{HZ}^\eta(W,Z;L)$ by taking into account, in \ref{c_HZ(W,Z;L)}, only Hamiltonian chords belonging to the class $\eta$. In the same way, one obtains an almost existence result for chords belonging to the class $\eta$ whenever this capacity is finite. We will use this version of the capacity for $\eta=0\in\pi_1(W,L)$.
\end{remark}

\begin{remark}
    As mentioned above, the absolute capacity introduced by Lisi and Rieser \cite[Definition 1.13]{c_coisotrope} is more general because it is defined for any coisotropic submanifold $N$ and any equivalence relation $\sim$ that is coarser than the one defined by the isotropic foliation. However, the definitions above can easily be adapted to this case, in order to define a relative coisotropic capacity $c_{HZ}(W,Z;N,\sim)$. Moreover, one can adapt the proofs of theorems \ref{c(W,Z;L) => cords} and \ref{c(W,Z;L) => presque existence} to show that the finiteness of this capacity will imply an almost-existence property for Hamiltonian chords with endpoints on the same leaf.
\end{remark}

\section{Conventions and overview of wrapped Floer homology}
\label{Convention and overview}

We give here an overview of wrapped Floer homology, based in large part on Ritter \cite{Ritter}.

\subsection{General setting}

We consider a Liouville domain $(W,\theta)$ and set $\omega=d\theta$. We consider an exact Lagrangian submanifold $L$ of $W$: $\theta_{|L}=df$, where $f$ is defined on $L$. We also assume that the intersection $\partial L = \partial W\cap L$ is transverse and that $\theta|_L=0$ in a neighborhood of $\partial L$. Let $Z$ denote the Liouville vector field on $W$ defined by $i_Z\omega = \theta$, and let $\alpha = \theta|_{\partial W}$ be the contact form induced on $\partial W$. The boundary $\partial L$ is then a Legendrian submanifold of the contact manifold $\partial W$. \\

The flow of $Z$ induces a parametrization $V\simeq (-\epsilon,0]\times\partial W$ of a neighborhood $V$ of $\partial W$. Writing $\rho$ for the coordinate on $(-\epsilon,0]$ and $r = e^\rho$, we may define the \emph{symplectic completion}
$$\hat{W} = W \cup [1,\infty)\times\partial W$$
where $r$ is the coordinate on $[1,\infty)$, and extend $Z,\theta$, and $\omega$ to $\hat{W}$ by
$$Z = \partial_\rho, \qquad \theta = r\alpha, \qquad \omega = d\theta$$
on $[1,\infty)\times\partial W$. Moreover, the Liouville vector field $Z$ is tangent to $L$ in the neighborhood of $\partial L$ where $\theta|_L=0$, so that near $\partial W$, $L$ takes the form $(1-\epsilon,1]\times\partial L$. We may then extend $L$ to a submanifold
$$\hat{L} = L\cup[1,\infty)\times\partial L$$
of the symplectic completion. We also extend the function $f$ to $\hat{L}$ by setting $f(r,x)=f(x)$ for $(r,x)\in[1,\infty)\times\partial L$. This function is locally constant on $\hat{L}\backslash L$ since $\theta|_{\hat{L}\backslash L}=0$. In particular, $\hat{L}$ is an exact Lagrangian submanifold of $\hat{W}$.\\

We will consider time-dependent Hamiltonians $H : \mathbb{R}\times W\to\mathbb{R}$, and write $H^t(x) = H(t,x)$. The Hamiltonian (time-dependent) vector field $X_{H_t}$ is defined by $\omega(X_{H_t},-) = dH^t$. Denoting by $\phi_t$ the flow of $H$, $H$ is said to be \emph{non-degenerate} if the Lagrangian submanifolds $\hat{L}$ and $\phi_{-1}(\hat{L})$ are transverse. A Hamiltonian $H$ is said to be \emph{linear at infinity} if, for $r$ sufficiently large,
$$\partial_rH = a$$
where $a$ is a constant called the \emph{slope at infinity}. It is said to be \emph{admissible for symplectic homology} if:
\begin{enumerate}
\item $H$ is linear at infinity (with slope $a$).
\item $\partial W$ admits no Reeb chords of period $a$ (or equivalently, $H$ admits no Hamiltonian chords in the region where it is linear).
\item $H$ is \emph{non-degenerate}.
\end{enumerate}
We denote by $\mathcal{H}$ the set of such Hamiltonians.\\

Let $(J^t)_{0\leq t\leq 1}$ be a family of almost complex structures on $\hat{W}$, calibrated by $\omega$. Such a family is said to be \emph{admissible} if, for $r$ sufficiently large, it is time-independent and cylindrical, that is:
\begin{enumerate}
    \item We have
    $$J^tZ=R,$$
    where $R$ is the Reeb vector field on $\partial W$.
    \item $J$ preserves the contact distribution $\text{Ker}(\theta|_{\partial W})$.
    \item The restriction of $J$ to this distribution is independent of $r$.
\end{enumerate}
We denote by $\mathcal{J}$ the set of such time-dependent complex structures.\\

These conditions provide the framework for Lagrangian Floer homology. More precisely, there exists a generic subset $(\mathcal{H}\times\mathcal{J})_\text{reg}\subset \mathcal{H}\times\mathcal{J}$ of \emph{regular pairs} such that Lagrangian Floer homology $FH(H,J)$ is defined for every regular pair $(H,J)$. \\

We will also need to vary the data $(H,J)$ through homotopies. Let $(H_\pm,J_\pm)$ be two admissible pairs. A homotopy from $(H_-,J_-)$ to $(H_+,J_+)$ is a smooth family $(H_s,J_s)_{s\in\mathbb{R}}$ such that $(H_s,J_s)=(H_-,J_-)$ for $s$ sufficiently close to $-\infty$, $(H_s,J_s)=(H_+,J_+)$ for $s$ sufficiently close to $+\infty$, and such that there exists $r_0\geq1$ with
$$\forall s\in\mathbb{R}, \forall r\geq r_0, \quad H_s(t,r,x) = a_sr+b_s \quad \text{and} \quad J_s^tZ = R$$
with $a_s,b_s\in\mathbb{R}$. The homotopy is said to be \emph{monotone at infinity} if $\partial_sa_s\leq0$. It is said to be \emph{monotone} if $\partial_sH_s\leq0$.\\

\subsection{Action functional}

If $M$ is a manifold and $N$ is a subset of $M$, we denote by $\Omega(M,N)$ the set of smooth maps $\gamma : [0,1] \to M$ such that $\gamma(0),\gamma(1)\in N$. To define Floer theory on $\hat{W}$, we work with the action functional:

\begin{align*}
    \mathcal{A}_H : \Omega(\hat{W},\hat{L}) & \longrightarrow \mathbb{R} \\
    \mathcal{A}_H(x) & = \int_{[0,1]}x^*\theta + f(x(0)) - f(x(1)) - \int_0^1H(t,x(t))dt.
\end{align*}

The critical points of this functional are exactly the Hamiltonian chords of $H$, of return-time 1. We denote by $\text{Crit}(\mathcal{A}_H)$ the set of such trajectories.\\

\begin{remark}
\label{Action trajectoires contractiles}
    Let $\Omega_0(\hat{W},\hat{L})$ be the set of trajectories with endpoints in $\hat{L}$ that are contractible relatively to $\hat{L}$. Let $x\in\Omega_0(\hat{W},\hat{L})$, and let
    \begin{align*}
        \tilde{x} : [0,1]^2 & \longrightarrow \hat{W} \\
        \tilde{x}(0,t) & \in \hat{L} \\
        \tilde{x}(1,t) & = x(t)
    \end{align*} 
    be a homotopy from a path in $\hat{L}$ to $x$. Now Stokes' theorem gives:
    $$\int_{[0,1]}x^*\theta+f(x(0))-f(x(1)) = \int_{[0,1]^2}\tilde{x}^*\omega = \int_0^1\int_0^1\omega(\partial_s\tilde{x},\partial_t\tilde{x})dsdt.$$
    Therefore, one useful expression of the action functional on contractible trajectories is given by
    \begin{align*}
    \mathcal{A}_H : \Omega_0(\hat{W},\hat{L}) & \longrightarrow \mathbb{R} \\
    \mathcal{A}_H(x) & = \int_{[0,1]^2}\tilde{x}^*\omega - \int_0^1H(t,x(t))dt.
\end{align*}
\end{remark}

\subsection{Spaces of trajectories}

To define the Floer trajectories, we use the $L^2$-metric given by $g(X,Y) = \int_0^1\omega(X(t),J^tY(t))dt$, and follow the \emph{positive} gradient of the action. This gives the following Floer equation:
$$\partial_su + J^t(u)\partial_tu + \nabla_{J^t}H^t(u) = 0$$
with boundary condition $u(s,0),u(s,1)\in \hat{L}$. The \emph{energy} of such a solution is given by
$$E(u) = \int_{\mathbb{R}\times[0,1]}|\partial_su|^2_{J^t}dtds \in [0,\infty]$$
and this energy is finite if and only if $u$ connects two chords $x_-$ and $x_+$. In this situation, we have
$$E(u) = \mathcal{A}_{H}(x_+) - \mathcal{A}_{H}(x_-),$$
where $x_\pm(t) = \underset{s\to\pm\infty}{\lim}u(s,t)$. We then define
$$\mathcal{M} = \{u \text{ solution of the Floer equation}, E(u)<\infty\} = \bigcup_{x_-,x_+\in\text{Crit}(\mathcal{A}_H)}\mathcal{M}(x_-,x_+)$$
where $\mathcal{M}(x_-,x_+)$ is the set of Floer solutions such that $u(s,t)\underset{s\to\pm\infty}{\longrightarrow}x_\pm(t)$.\\

In the framework that we consider, it is not possible to directly define the index of a critical point of the action. However, we can associate an integer to each trajectory $u\in\mathcal{M}(x_-,x_+)$ (playing the role of the index difference between $x_-$ and $x_+$), this integer possibly depending on the trajectory in question. This construction is due to Viterbo \cite[Definition 1]{Indice_Maslov_Viterbo}. As a consequence, the resulting Floer homology will not be graded, but this property has no impact on our applications.\\

Let $x_\pm$ be two chords. We define the space of Floer solutions of index $k$:
$$\mathcal{M}_k(x_-,x_+) = \{u\in\mathcal{M}(x_-,x_+), m(u) = k\}.$$

If $(H,J)$ is a regular pair, then $\mathcal{M}_k(x_-,x_+)$ is a smooth manifold of dimension $k$. We also define the spaces of Floer trajectories (or unparametrized solutions):
\begin{align*}
\mathcal{L}(x_-,x_+) =& \mathcal{M}(x_-,x_+)/\mathbb{R}\\
\mathcal{L}_k(x_-,x_+) =& \mathcal{M}_{k+1}(x_-,x_+)/\mathbb{R}
\end{align*}
where $\mathbb{R}$ acts on $\mathcal{M}$ via
$$(r\cdot u)(s,t) = u(s+r,t).$$

For all $k\geq0$, $\mathbb{R}$ acts freely on $\mathcal{M}_{k+1}(x_-,x_+)$, and thus $\mathcal{L}_k(x_-,x_+)$ is a manifold of dimension $k$. These spaces can be compactified into \emph{manifolds with corners} (see \cite[Definition A.3]{Barraud_Cornea}) $\overline{\mathcal{L}}_{k}(x_-,x_+)$, by adding broken trajectories (see \cite[Appendix A]{Barraud_Cornea} for a proof):
$$\partial\overline{\mathcal{L}}_{k-1}(x_-,x_+) = \bigcup_{\underset{k_0+...+k_r=k}{0\leq k_0,...,k_r<k}}\bigcup_{y_1,...,y_r}\mathcal{L}_{k_0-1}(x_-,y_1)\times\mathcal{L}_{k_1-1}(y_1,y_2)\times...\times\mathcal{L}_{k_r-1}(y_r,x_+).$$

\subsection{Floer homology}

The Floer homology $FH_*(H,J)$ associated to a regular pair $(H,J)$ is the homology of the $\mathbb{Z}/2\mathbb{Z}$-module $FC_*(H,J)$ generated by the chords of $H$, equipped with the differential:

\begin{align*}
d : FC_*(H,J) & \longrightarrow FC_*(H,J) \\
d\langle x \rangle & = \sum_{x_-\in\text{Crit}(\mathcal{A}_H)}\#\mathcal{L}_0(x_-,x)\langle x_-\rangle.
\end{align*}

We can define a submodule $FC^{<a}_*(H,J)\subset FC_*(H,J)$ by restricting to chords of action strictly less than $a$, and one checks that the differential $d$ preserves this submodule. We also define, for $a\leq b$,
$$FC^{[a,b)}_*(H,J) = FC^{<b}_*(H,J)/FC^{<a}_*(H,J).$$

The homology of $FC^{<a}_*(H,J)$ (resp. $FC^{[a,b)}_*(H,J)$) will be denoted $FH^{<a}_*(H,J)$ (resp. $FH^{[a,b)}_*(H,J)$).\\

\subsection{Wrapped Floer homology}

The Lagrangian symplectic homology (or wrapped Floer homology) is constructed as a colimit of Floer homologies, using continuation maps that come from homotopies: each homotopy $(H_s,J_s)$ that is monotone at infinity induces a continuation map
\begin{align*}
    \psi : FC_*(H_+,J_+) & \longrightarrow FC_*(H_-,J_-)\\
    \psi\langle x_+\rangle & = \sum_{x_-\in\text{Crit}(\mathcal{A}_{H_-})}\#\mathcal{M}^{(H_s,J_s)}_0(x_-,x_+)\langle x_-\rangle.
\end{align*}
where $\mathcal{M}^{(H_s,J_s)}_k(x_-,x_+)$ are the spaces of solutions of the parametrized Floer equation
$$J_s^t(u)\partial_tu + \nabla_{J_s^t}H_s^t(u) + \partial_su = 0.$$

Observe that the Floer homology $FH_*(H,J)$ only depends on the slope of $H$ at infinity: if $H_-$ and $H_+$ have the same slope at infinity, then the continuation map $\psi$ is invertible, and its inverse is the continuation map associated to the time-reversed homotopy. Furthermore, we have an isomorphism
$$\Psi_H : FH_*(H,J) \underset{\simeq}{\longrightarrow} SH^{<a}_*(W,L)$$
where $a$ is the slope of $H$ at infinity. This allows to write wrapped homology as a colimit on the slopes:
$$SH(W,L) = \underset{a}{\text{colim }}SH^{<a}(W,L).$$

In part \ref{Le critère homologique}, we will focus on contractible chords (our criterion is only able to detect these chords). Therefore, for each homotopy cass $\eta\in\pi_1(W,L)\simeq\pi_1(\hat{W},\hat{L})$, we define Floer homology $FH_{*,\eta}(H,J)$ and symplectic homology $SH_{*,\eta}(W,L)$ by the same constructions, taking only into account the Hamiltonian chords that belong to the class $\eta$, and we have the decomposition:
$$FH_*(H,J) = \bigoplus_{\eta\in\pi_1(W,L)}FH_{*,\eta}(H,J).$$
The same holds for (filtered) wrapped Floer homology.

\subsection{From symplectic homology to Morse homology}
\label{Morse -> symplectique}

Let us consider a Weinstein neighbourhood of $\hat{L}$ in $\hat{W}$, which can be symplectically identified with the cotangent bundle $T^*\hat{L}$. It is a common fact that for $H$ sufficiently $\mathcal{C}^2$-small, the displaced Lagrangian $\phi_{-1}(\hat{L})$ can be identified with the graph of an exact 1-form $df_H$ in $T^*\hat{L}$ (see \cite[Proposition 9.4.2]{McDuffSalamon17}). Moreover, since the intersection is transverse, $f_H$ is a Morse function and its critical points correspond to the intersection points of $\hat{L}$ and $\phi_{-1}(\hat{L})$. According to \cite[Proposition 4.6]{Oh}, given a metric $g$ such that $(-f_H,-\nabla_gf_H)$ is Morse-Smale, there exists an almost complex structure $J_{g,H}$ such that $(H,J_{g,H})$ is regular, and we have an isomorphism
$$FH_*(H,J_{g,H})\simeq MH_*(\hat{L};-f_H,-\nabla_gf_H)\simeq H_*(L,\partial L),$$
where $MH_*$ stands for the Morse homology. The original result is due to Floer \cite{Floer_Wittenscomplex}. Since the left-hand side is a part of the colimit construction, we get a canonical morphism $i_L : H_*(L,\partial L)\to SH_*(W,L)$, which is independent of $H$ and factors through all filtered homology groups $SH^{<a}_*(W,L)$:

\[\begin{tikzcd}
    & H_*(L,\partial L) \arrow[r,"i_a"]\arrow[rd,"i_L"] & SH_*^{<a}(W,L) \arrow[d,"k_a"] \\
    & & SH_*(W,L)
\end{tikzcd}\]

\begin{remark}
\label{i_L factors through SH0}
    Since the Hamiltonian chords that are contained in the Weinstein neighborhood of $L$ are contractible relatively to $L$, we observe that the morphism $i_L$ (resp. $i_a$) also factors through $SH_{*,0}(W,L)\hookrightarrow SH_*(W,L)$ (resp. $SH_{*,0}^{<a}(W,L)\hookrightarrow SH_*^{<a}(W,L)$).
\end{remark}

\section{The homological criterion}
\label{Le critère homologique}

This section is devoted to the proof of Theorem \ref{critère homologique => cordes}, which provides a homological criterion for the almost existence of Hamiltonian chords near a hypersurface. This is the Lagrangian version (with non-DG coefficients) of Theorem 5.1 in \cite{dg-floer2}. This criterion involves the morphism $i_L$ defined earlier between Morse homology of the Lagrangian submanifold $L$ and wrapped Floer homology, as well as the shriek morphism in Morse homology, summarized in the diagram below:

\[\begin{tikzcd}
    & H_*(L,\partial L) \arrow[r,"i_L"]\arrow[d,"j_{U!}"] & SH_{*,0}(W,L) \\
    & H_*(U\cap L,\partial U\cap L) &
\end{tikzcd}\]

Here, $U$ is a relatively compact open subset of $W$, and $j_U$ is the inclusion $U\cap L\hookrightarrow L$. We recall that Theorem \ref{critère homologique => cordes} asserts that the almost-existence property holds for each hypersurface of the form $\Sigma=\partial U$ for which the homological criterion $\text{Ker}(i_L)\not\subset\text{Ker}(j_{U!})$ is satisfied. To establish this result, we use the relative Lagrangian Hofer–Zehnder capacity defined in Section \ref{Capacités et orbites}, which we will bound by a capacity defined from wrapped Floer homology, following the ideas of \cite[Section 5.1]{dg-floer2} and \cite[Section 3.3]{BK22}.

\subsection{A spectral invariant}

 We begin by adapting the spectral invariant defined in \cite{dg-floer2} to the Lagrangian case. Let $H : \mathbb{R}\times\hat{W}\to\mathbb{R}$ be an admissible Hamiltonian of slope $a$ at infinity, and $\sigma\in SH^{<a}_{*,0}(W,L)$. We define:
$$\rho(H,\sigma) = \inf\{b, \Psi^{-1}_H(\sigma)\in\text{Im}(FH^{<b}_{*,0}(H)\to FH_{*,0}(H))\}$$
where $\Psi_H$ is the isomorphism $FH_{*,0}(H)\simeq SH_{*,0}^{<a}(W,L)$. The basic properties of this invariant are analogous to the case of periodic orbits, and are summarized in the following proposition.

\begin{proposition}
\label{propriétés invariant spectral}
Let $\sigma\in SH^{<a}_{*,0}(W)$. Then:
\begin{enumerate}
\item (Continuity) The function $\rho(-,\sigma)$ extends to a continuous function on the space of Hamiltonians (not necessarily admissible) of slope $a$ at infinity, equipped with the Hofer norm:
$$||H|| = \int_0^1\left(\underset{\hat{W}}{\max}H^t-\underset{\hat{W}}{\min}H^t\right)dt.$$
\item (Spectrality) For any Hamiltonian $H$ (not necessarily admissible) of slope $a$ at infinity,
$$\rho(H,\sigma)\in\text{\emph{Spec}}_0(H):=\{\mathcal{A}_H(x), x\text{ contractible Hamiltonian chord}\}.$$
\end{enumerate}
\end{proposition}

\textit{Proof}.

\textit{Continuity.} Let $H_-,H_+$ be two admissible Hamiltonians of slope $a$ at infinity. Consider a homotopy $(H_s,J_s)$ that is monotone at infinity and such that $H_s\underset{s\to\pm\infty}{\to}H_\pm$. Let $u : \mathbb{R}\times[0,1]\to \hat{W}$ be a solution of the parametrized Floer equation:
$$\partial_s u + J_s^t(u)\partial_tu + \nabla_{J_s^t}H_s^t(u) = 0$$
$$u(s,0),u(s,1)\in L$$
connecting two trajectories $x_-$ and $x_+$ (Hamiltonian chords of $H_-$ and $H_+$, respectively). Let $R=\overline{\mathbb{R}}\times[0,1]$. We have:
\begin{align*}
\int_R|\partial_su|^2_{J_s}dtds & = \int_R\langle\partial_su,-J_s^t(u)\partial_tu - \nabla_{J_s^t}H_s^t(u)\rangle dtds\\
& = \int_R\omega(\partial_su,\partial_tu)dtds - \int_RdH_s^t(\partial_su)dtds\\
& = -\int_Ru^*\omega - \int_R\frac{d}{ds}H_s^t(u)dsdt + \int_R\partial_sH_s^t(u)dtds\\
& = -\int_Ru^*\omega + \int_0^1[H_-^t(x_-(t))-H_+^t(x_+(t))]dt + \int_R\partial_sH_s^t(u)dtds.
\end{align*}

Since the symplectic form is assumed to be exact ($\omega = d\theta$), we obtain from Stokes’ theorem
\begin{align*}
\int_Ru^*\omega & = \int_{\partial R}u^*\theta\\
& = \int_{[0,1]}x_-^*\theta + \int_{\overline{\mathbb{R}}}u(-,1)^*\theta - \int_{[0,1]}x_+^*\theta - \int_{\overline{\mathbb{R}}}u(-,0)^*\theta.
\end{align*}

Now note that $u(-,1)$ and $u(-,0)$ take values in $\hat{L}$, and that $\theta = df$ on $\hat{L}$, which allows us to conclude:
\begin{align*}
    \int_Ru^*\omega & = \int_{[0,1]}x_-^*\theta - \int_{[0,1]}x_+^*\theta + f(x_+(1)) - f(x_+(0)) - f(x_-(1)) + f(x_-(0))\\
    & = \mathcal{A}_{H_-}(x_-) - \mathcal{A}_{H_+}(x_+) + \int_0^1[H_-^t(x_-(t))-H_+^t(x_+(t))]dt.
\end{align*}

Thus, we obtain the following identity:
\begin{equation}
\label{identité homotopie d'hamiltoniens}
0\leq \int_\mathbb{R}\int_{[0,1]}|\partial_su|^2_{J_s}dtds = \mathcal{A}_{H+}(x_+) - \mathcal{A}_{H-}(x_-) + \int_{\mathbb{R}}\int_{[0,1]}\partial_sH_s^t(u)dtds.
\end{equation}

Moreover, since the non-degeneracy property required for the homotopy $(H_s,J_s)$ is generic, we can choose a homotopy arbitrarily close\footnote{For the needs of the proof, one perturbs $\overline{H}_s$ into a homotopy $H_s = \overline{H}_s + h_s$ with
$$\int_\mathbb{R}\int_0^1\underset{\hat{W}}{\max}(\partial_sh_s^t)dtds \leq\epsilon$$
for arbitrarily small $\epsilon$.} to
$$\overline{H}_s = H_- + \beta(s)(H_+-H_-),$$
where $\beta : \mathbb{R}\to[0,1]$ is a non-decreasing function that is equal to 0 near $-\infty$ and equal to 1 near $+\infty$. We then obtain an upper bound of the form
$$\int_{\mathbb{R}}\int_{[0,1]}\partial_sH_s^t(u)dtds \leq \int_0^1\underset{\hat{W}}{\max}(H_+^t-H_-^t)dt+\epsilon=: E,$$
and therefore
$$\mathcal{A}_{H_-}(x_-)\leq\mathcal{A}_{H_+}(x_+)+E.$$

This shows that the corresponding continuation map $\psi$ sends the subcomplex $FC^{<a}_{*,0}(H_+)$ to $FC_{*,0}^{<a+E}(H_-)$. We thus have the following commutative diagram:
\[\begin{tikzcd}
    & FH^{<b}_{*,0}(H_+) \arrow[d,"\psi"]\arrow[r] & FH_{*,0}(H_+) \arrow[r,"\Psi_{H_+}"]\arrow[d,"\psi"] & SH_{*,0}^{<a}(W,L) \arrow[ld,leftarrow,"\Psi_{H_-}"] \\
    & FH_{*,0}^{<b+E}(H_-) \arrow[r] & FH_{*,0}(H_-) &
\end{tikzcd}\]
where the three morphisms in the right triangle are isomorphisms. Since $\epsilon$ can be chosen arbitrarily small, we obtain
$$\rho(H_-,\sigma)-\rho(H_+,\sigma) \leq \int_0^1\underset{\hat{W}}{\max}(H_+^t-H_-^t)dt,$$
and by exchanging the roles of $H_+$ and $H_-$, we finally deduce
$$|\rho(H_-,\sigma)-\rho(H_+,\sigma)|\leq ||H_--H_+||,$$
where $||-||$ is the Hofer norm. Hence the continuity property of the spectral invariant.\\

\textit{Spectrality.} Assume that $H$ is non-degenerate. Since $\text{Spec}_0(H)$ is closed, if $b\notin\text{Spec}_0(H)$, there exists a neighborhood $(b-\epsilon,b+\epsilon)$ of $b$ contained in $\mathbb{R}\backslash\text{Spec}_0(H)$. For any value $c$ in this neighborhood, we have $FH_{*,0}^{<c}(H)\simeq FH_{*,0}^{<b}(H)$, which shows that $b$ cannot be the infimum defining the spectral invariant.\\

Now consider any Hamiltonian $H$, which we can approximate by a sequence of non-degenerate Hamiltonians $(H_k)$ in the $\mathcal{C}^2$-norm:
$$||H||_{\mathcal{C}^2} = ||H||_\infty + \underset{t\in\mathbb{R}}{\sup}||X_H||_\infty + \underset{||X||_\infty\leq1,t\in\mathbb{R}}{\sup}||\nabla_XX_H||_\infty$$
where vector fields are measured with respect to an arbitrary Riemannian metric $g$ on $\hat{W}$, and $\nabla$ is the associated Levi-Civita connection. Moreover, outside a compact set $K$, $H$ is linear and we can choose the sequence $(H_k)$ such that $H_k=H$ outside $K$. For each $k$, there exists a Hamiltonian chord $x_k$ such that $\rho(H_k,\sigma)=\mathcal{A}_{H_k}(x_k)$. Furthermore, since the trajectories $x_k$ are contained in the compact set $K$, the estimate
\begin{align*}
d(x_k(t_1),x_k(t_0)) & \leq \int_{t_0}^{t_1}|\dot{x}_k(t)|dt\\
& = \int_{t_0}^{t_1}|X_{H_k^t}(x_k(t))|dt\\
& \leq ||H_k||_{\mathcal{C}^2}|t_1-t_0|
\end{align*}
allows us to apply Arzelà-Ascoli theorem to the family $(x_k)$. There exists therefore a subsequence of $(x_k)$ converging to a trajectory $x$ in the $\mathcal{C}^0$ sense. Moreover, the family $(\dot{x}_k)$ is bounded, and denoting by $P_x^y$ the parallel transport from $x$ to $y$, we have:
\begin{align*}
    d(\dot{x}_k(t_0),\dot{x}_k(t_1)) := & |P_{x_k(t_0)}^{x_k(t_1)}(\dot{x}_k(t_0))-\dot{x}_k(t_1)|\\
    = & \left| \int_{t_0}^ {t_1} P_{x_k(t)}^{x_k(t_1)}(\nabla_{\dot{x}_k(t)}\dot{x}_k(t))dt \right| \\
    \leq & \int_{t_0}^ {t_1} |\nabla_{X_{H_k^t}(x_k(t))}X_{H_k^t}(x_k(t))|dt \\
    \leq & (t_1-t_0)\underset{t\in\mathbb{R}}{\sup}||X_{H_k^t}||_\infty\underset{||X||\leq\infty,t\in\mathbb{R}}{\sup}||\nabla_XX_{H_k^t}||_\infty \\
    \leq & (t_1-t_0)||H_k||_{\mathcal{C}^2}^2.
\end{align*}

Hence $(\dot{x}_k)$ is also equicontinuous, and the Arzelà–Ascoli theorem shows that $x$ is of class $\mathcal{C}^1$ and that $x_k$ converges to $x$ in $\mathcal{C}^1$-norm. Thus, it suffices to take the limit in the equation
$$\dot{x}_k = X_{H_k}(x_k)$$
to deduce that $x$ is a trajectory of $H$. Since the Lagrangian $\hat{L}$ is closed in $\hat{W}$, the trajectory $x$ moreover satisfies $x(0),x(1)\in L$. By continuity of the spectral invariant, we have:
$$\rho(H,\sigma) = \mathcal{A}_H(x).$$

\qed

\subsection{The spectral capacity}

Let us consider a compactly supported Hamiltonian $H : \mathbb{R}\times W\to\mathbb{R}$. Fix $\delta>0$ and $r_0>1$. We can extend $H$ to a smooth function $H_\delta : \hat{W}\to\mathbb{R}$, with:
$$H_\delta(r,x) = h_\delta(r), \forall r\geq1$$
where $r$ is the coordinate on $[1,\infty)$ in $\hat{W}=W\cup[1,\infty)\times\partial W$, and $h_\delta : [1,\infty)\to\mathbb{R}$ is non-decreasing, strictly convex on $(1,r_0)$ and linear of slope $\delta$ on $[r_0,\infty)$.\\

Recall that for every $a>0$, the morphism $i_L$ between Morse homology of $L$ and symplectic homology factors as
\[\begin{tikzcd}
    & H_*(L,\partial L) \arrow[r,"i_a"]\arrow[rd,"i_L"] & SH^{<a}_{*,0}(W,L) \arrow[d,"k_a"] \\
    & & SH_{*,0}(W,L)
\end{tikzcd}\]
where $k_a$ is the natural morphism
$$SH^{<a}_{*,0}(W,L) \longrightarrow \underset{b}{\text{colim }}SH_{*,0}^{<b}(W,L)\simeq SH_{*,0}(W,L).$$
Now if $\alpha\in H_*(L,\partial L)$ is a Morse homology class, we define
$$\sigma(H,\alpha) := \underset{\delta\to0}{\lim}\rho(H_\delta,i_\delta(\alpha)).$$

\begin{definition}
    The \emph{absolute} and \emph{relative spectral capacity} are respectively given by
    $$c_\sigma(W;L) := \underset{H\in C^\infty_c(\mathbb{R}\times W)}{\sup}\sigma(H,[L,\partial L]),$$
    $$c_\sigma(W,U;L) := \underset{H\in C^\infty_c(\mathbb{R}\times W)}{\sup}\inf\{\sigma(H,\alpha), \alpha\in H_*(L,\partial L)\backslash\text{Ker}(j_{U!})\}.$$
\end{definition}

Recall that the Hofer-Zehnder capacity was defined in part \ref{Capacités et orbites} using a variational principle on admissible Hamiltonians $H\in\mathcal{H}(W,Z;L)$, and that we denote by $\underline{\mathcal{H}}(W,Z;L)$ the set of such Hamiltonians that do not have non-constant chords of period less than 1. One key result in the proof of Theorem \ref{critère homologique => cordes} is the following:

\begin{proposition}
\label{rho(H_delta)}
    Let $H\in\underline{\mathcal{H}}(W,\overline{U};L)$. Then for $\delta$ sufficiently small and for each class $\alpha$ such that $j_{U!}(\alpha)\neq0$,
    $$\rho(H_\delta,i_\delta(\alpha)) = -\min(H).$$
\end{proposition}

An immediate corollary of this proposition is the bounding of the Hofer-Zehnder capacity by the spectral capacity:

\begin{corollary}
\label{c_HZ < c_s}
    For any relatively compact open set $U$ of $W$,
    $$c_{HZ}^0(W,\overline{U};L)\leq c_\sigma(W,U;L)$$
\end{corollary}

Here $c_{HZ}^0$ is the Hofer-Zehnder capacity that only takes into accounts the contractible chords, as mentioned in Remark \ref{pi_1 sensitive capacities}, also called $\pi_1$-sensitive Hofer-Zehnder capacity.\\

The rest of this section is devoted to the proof of Proposition \ref{rho(H_delta)}. To this end we will use the continuity and spectrality properties of the invariant $\rho$ to reduce to the computation of $\rho$ for a rescaled Hamiltonian $\epsilon H_{\delta/\epsilon}$, for which the Floer complex coincides with a Morse complex. However, a difference appears compared with the non-Lagrangian case: while it is true that, in the Lagrangian setting, a sufficiently small Hamiltonian $H$ induces a Floer complex isomorphic to a Morse complex, the Morse function in question is not directly $H$, and not all chords need to be constant (see Section \ref{Morse -> symplectique}). We therefore use the following lemma:

\begin{lemma}
\label{Estimée métrique}
    Let $L$ be an exact Lagrangian submanifold of a Liouville domain $(W,\theta)$. Fix a Riemannian metric $g$ on $W$. There exists a constant $C>0$ such that, for a sufficiently $\mathcal{C}^1$-small Hamiltonian $H$ on $W$, for any 1-periodic Hamiltonian chord $x$,
    $$\left|\mathcal{A}_H(x) + \int_0^1H(t,x(t))dt\right|\leq C\underset{t\in[0,1]}{\sup}||X_{H^t}||_\infty^2.$$
\end{lemma}

\textit{Proof.} Let us cover the compact manifold $L$ by a finite collection of open balls $B_i = B(x_i,r_i)$ such that $x_i\in L$ and $\overline{B_i}$ is contained in the source of a Darboux chart:
$$\Phi_i : (U_i,U_i\cap L) \longrightarrow (\mathbb{R}^{2n},\mathbb{R}^n)=(\{q^j,p_j\},\{q^j\})$$
(or $\Phi_i : (U_i,U_i\cap L) \longrightarrow (\{q^j,p_j\},\{q^j\})\cap\{q^1\leq0\}$ if $x_i\in\partial L$).\\

For $H$ sufficiently $\mathcal{C}^1$-small, any 1-periodic chord is contained in one of the open balls $B_i$, and is, in particular, contractible. Therefore, we use the expression of the action functional given in Remark \ref{Action trajectoires contractiles}. Now let us specify the homotopies we use. For each Hamiltonian chord $x$, there exists $i$ such that $x$ is contained in $B_i$, so $\Phi_i\circ x$ is a chord of $K_i := H\circ\Phi_i^{-1}$ on $(\mathbb{R}^{2n},\mathbb{R}^{n})$. Writing $\Phi_i\circ x(t) = y(t) = (q(t),p(t))$, consider the following homotopy:
\begin{align*}
    \tilde{y} : [0,1]^2 & \longrightarrow \mathbb{R}^{2n}\\
    (s,t) & \mapsto q(t) + sp(t).
\end{align*}
Then the homotopy we use for the chord $x$ is given by
$$\tilde{x} := \Phi^{-1}_i\circ\tilde{y}.$$
The quantity we want to estimate is
\begin{align*}
    \mathcal{A}_H(x)+\int_0^1H(t,x(t))dt & = \int_{[0,1]^2}\tilde{x}^*\omega \\
    & = \int_{[0,1]^2}\tilde{y}^*(\Phi_i^{-1})^*\omega \\
    & = \int_{[0,1]^2}\tilde{y}^*\omega_0 \\
    & = \int_0^1\int_0^1\omega_0(\partial_s\tilde{y},\partial_t\tilde{y})dtds,
\end{align*}
where $\omega_0$ is the standard symplectic form on $\mathbb{R}^{2n}$. Let $J_0$ denote the canonical (almost) complex structure on $\mathbb{R}^{2n}$. Then
$$|\omega_0(\partial_s\tilde{y},\partial_t\tilde{y})| = |\langle J_0\partial_s\tilde{y},\partial_t\tilde{y}\rangle|\leq |J_0\partial_s\tilde{y}|\cdot|\partial_t\tilde{y}| = |\partial_s\tilde{y}|\cdot|\partial_t\tilde{y}|$$
(here the norms must be understood in the sense of the standard metric on $\mathbb{R}^{2n}$, that is $\omega_0(J_0-,-)$), hence the estimate
$$\left|\mathcal{A}_H(x)+\int_0^1H(t,x(t))dt\right|\leq ||\partial_s\tilde{y}||_\infty||\partial_t\tilde{y}||_\infty.$$
Now,
\begin{align*}
    |\partial_s\tilde{y}| & = |p(t)| \leq\int_0^t|\dot{p}(t)|dt\leq ||X_K||_\infty,\\
    |\partial_t\tilde{y}| & = |\dot{q}(t)+s\dot{p}(t)| \leq |\dot{y}(t)|\leq ||X_K||_\infty,
\end{align*}
so
$$\left|\mathcal{A}_H(x)+\int_0^1H(t,x(t))dt\right| \leq ||X_K||_\infty^2 = \underset{t\in[0,1]}{\sup}\underset{B_i}{\sup}||X_{H^t}||^2_{\Phi_i^*g_0} \leq C_i\underset{t\in[0,1]}{\sup}||X_{H^t}||^2_{g,\infty}$$
where the constant $C_i$, depending only on $i$, comes from the fact that all the metrics are equivalent in the compact set $\overline{B_i}$. Since there are only a finite number of balls $B_i$, we finally get the desired estimate.

\qed

Now let us consider the rescaled Hamiltonians $\epsilon H_{\delta/\epsilon}$, for $0<\epsilon\leq1$. We fix a metric $g$ on $\hat{W}$.

\begin{lemma}
\label{estimation rho(epsilon H)}
    There exists a constant $C$ such that, for all $\delta,\epsilon$ small enough, we have the estimate :
    $$\left|\rho(\epsilon H_{\delta/\epsilon},i_\delta(\alpha)) + \epsilon\min(H)\right| \leq \epsilon^2C||X_H||^2_g.$$
\end{lemma}

\textit{Proof}. Let $\delta>0$ be smaller than all the periods of Reeb chords on $(\partial W,\partial L)$, and let $\epsilon>0$ be such that $2\epsilon H$ is $\mathcal{C}^1$-small. Then $\epsilon H_{\delta/\epsilon}$ can be approximated by a sequence of Hamiltonians $H_k\in\mathcal{H}, k\geq0$, such that for every integer $k$, $H_k$ is $\mathcal{C}^1$-small and equal to $H$ in the region ${r\geq r_k}$, with $1\geq r_k\to1$ as $k\to\infty$. By Lemma \ref{Estimée métrique} applied to the truncated manifold $\hat{W}(r_k) = W\cup\partial W\times[1,r_k]$, there exists a constant $C$ such that, after possibly reducing $\epsilon$ so that $\epsilon H_{\delta/\epsilon}$ and the $H_k$ (for $k$ large enough) are sufficiently small, we have, for large $k$ and for every chord $x\in\text{Crit}(\mathcal{A}_{H_k})$,
\begin{equation}
\label{estimation A_H_k}
    \left|\mathcal{A}_{H_k}(x)+\int_0^1 H_k(x(t))dt\right|\leq C\underset{\hat{W}(r_k)}{\sup}|X_{H_k}|^2\underset{k\to\infty}{\longrightarrow}\epsilon^2C||X_H||^2.
\end{equation}

Recall that, since the Hamiltonians $H_k$ are $\mathcal{C}^1$-small, there exist Morse functions $f_k$ on $\hat{L}$ such that every critical point $y$ of $f_k$ is a point of the intersection $L\cap\phi_{-1}^{H_k}(L)$, corresponding to the chord $x(t) = \phi_t^{H_k}(y) =: \Psi_k(y)(t)$, and this correspondence is bijective. Choose the Hamiltonians $H_k$ so that $\nabla f_k$ points inward along $\partial\overline{U}$. By definition of the shriek morphism, the assumption $j_{U!}(\alpha)\neq0$ means precisely that $\alpha$ is not entirely represented by critical points outside $U$ (this being valid for all the Morse functions $f_k$). By definition of the spectral invariant, this means that:
$$\rho(H_k,i_\delta(\alpha))\geq\min(\mathcal{A}_H\circ\Psi_k(\text{Crit}(f_k)\cap U)).$$

Since $H$ is admissible for $c_{HZ}^0(W,\overline{U};L)$, we have $H=\min(H)$ on a neighborhood $V$ of $\overline{U}$. For large $k$, we may then assume that, for every point $y\in U$, the corresponding trajectory $x=\Psi_k(y)$ is entirely contained in $V$. If moreover $y$ is a critical point of $f_k$, the estimate (\ref{estimation A_H_k}) remains valid when the $\mathcal{C}^0$-norm of $X_{H_k}$ is evaluated only on $V$:
$$\left|\mathcal{A}_{H_k}(x)+\int_0^1 H_k(x(t))dt\right|\leq C\underset{U}{\sup}|X_{H_k}|^2 \underset{k\to\infty}{\longrightarrow}0.$$

Since $H_k\to\epsilon\min(H)$ on $U$, we deduce that
$$\mathcal{A}_H\circ\Psi_k(\text{Crit}(f_k)\cap U) \subset [-\epsilon\min(H)-\eta_k,-\epsilon\min(H)+\eta_k]$$
with $\eta_k\to0$ as $k\to\infty$. In particular, $\rho(H_k,i_\delta(\alpha))\geq-\epsilon\min(H)-\eta_k$, and therefore, taking the limit $k\to\infty$:
\begin{equation}
    \label{minoration rho}
    \rho(\epsilon H_{\delta/\epsilon},i_\delta(\alpha))\geq-\epsilon\min(H).
\end{equation}

On the other hand, we can use estimate (\ref{estimation A_H_k}) to bound the action of the set of chords of $H_k$:
$$\forall x\in\text{Crit}(\mathcal{A}_H), \mathcal{A}_H(x)\leq -\min(H_k) + C\underset{\hat{W}(r_k)}{\sup}||X_{H_k}||^2,$$
so that
$$\rho(H_k,i_\delta(\alpha))\leq -\min(H_k) + C\underset{\hat{W}(r_k)}{\sup}||X_{H_k}||^2.$$

Taking the limit $k\to\infty$, we obtain:
\begin{equation}
    \label{majoration rho}
    \rho(\epsilon H_{\delta/\epsilon},i_\delta(\alpha))\leq -\epsilon\min(H) + \epsilon^2C||X_H||^2.
\end{equation}

Combining (\ref{minoration rho}) and (\ref{majoration rho}) finally yields the desired estimate.

\qed

\textit{Proof of Proposition \ref{rho(H_delta)}}. Recall that, by assumption, $H$ admits no non-constant contractible chord of return time less than $1$. Since $\epsilon H_{\delta/\epsilon}$ likewise admits no non-constant chord (of return time $1$) in $\hat{W}\setminus W$ for every $0<\epsilon$, we deduce that
$$\forall 0<\epsilon\leq1, \text{Spec}_0(\mathcal{A}_{\epsilon H_{\delta/\epsilon}}) = \epsilon\text{Spec}_0(\mathcal{A}_{H_\delta}) =: \epsilon\mathcal{E},$$
where the spectrum $\mathcal{E}$ consists only of critical values of $H_\delta$. The function $\epsilon\mapsto\frac{1}{\epsilon}\rho(\epsilon H_{\delta/\epsilon},i_\delta(\alpha))$ is continuous (by continuity of $\rho$), and takes values in $\mathcal{E}$ (by spectrality of $\rho$). Now, $\mathcal{E}$ is totally disconnected by Sard’s theorem, hence this function is constant. Thus,
$$\forall\epsilon\in(0,1], \rho(\epsilon H_{\delta/\epsilon},i_\delta(\alpha)) = \epsilon\rho(H_\delta,i_\delta(\alpha)).$$

Lemma \ref{estimation rho(epsilon H)} then yields, for $\epsilon$ sufficiently small:
$$\left|\rho(H_\delta,i_\delta(\alpha)) + \min(H)\right|\leq \epsilon C||X_H||^2_g,$$
and therefore
$$\rho(H_\delta,i_\delta(\alpha)) = -\min(H).$$

\qed
\\

\subsection{The capacity $c_{SH}$ and the homological criterion}

Let $\alpha\in H_*(L,\partial L)$ be a Morse homology class. We define
$$c(\alpha) = \inf\{a\geq0, i_a(\alpha)=0\} \in [0,+\infty].$$
In particular, $c(\alpha)<\infty$ if and only if $i_L(\alpha)=0$.

\begin{definition}
\label{SH capacity}
    We define the numbers
    $$c_{SH}(W;L) := c([L,\partial L]) = \sup\{c(\alpha), \alpha\in H_*(L,\partial L)\},$$
    $$c_{SH}(W,U;L) := \inf\{c(\alpha), \alpha\in H_*(L,\partial L)\backslash\text{Ker}(j_{U!})\}.$$
\end{definition}

The second equality in the definition of $c_{SH}(W;L)$ comes from the fact that the maps $i_a$ are ring morphisms: for all $a,\epsilon>0$ and all classes $\alpha\in H_*(L,\partial L)$, we have
    $$i_{a+\epsilon}(\alpha) = i_{a+\epsilon}(\alpha\cdot[L,\partial L]) = i_\epsilon(\alpha)\cdot i_a([L,\partial L]),$$
so if $i_a([L,\partial L])=0$ then $i_{a+\epsilon}=0$ for arbitrary small $\epsilon>0$.\\

An important observation is that $c_{SH}(W;L)<\infty$ if and only if the wrapped Floer homology $SH_*(W,L)$ vanishes, and $c_{SH}(W,U;L)<\infty$ if and only if there exists a class $\alpha$ satisfying the homological criterion described in Theorem \ref{critère homologique => cordes}.

\begin{proposition}
\label{c_s < c_SH}
    For each open set $U$ of $W$,
    $$c_\sigma(W,U;L)\leq c_{SH}(W,U;L).$$
\end{proposition}

\textit{Proof}. We will prove that for any $H\in C^\infty_c(\mathbb{R}\times W)$ and any class $\alpha\notin\text{Ker}(j_{U!})$, we have
$\sigma(H,\alpha)\leq c(\alpha).$
Reason by contradiction and assume that there exists such a pair $(H,\alpha)$ with $\sigma(H,\alpha)> c(\alpha)$. Then we can find real numbers $a,\delta$ such that
$$\rho(H_\delta,i_\delta(\alpha))>a>c(\alpha)>\delta$$
(it suffices to take $\delta$ sufficiently small). We consider perturbations (in the sense of the Hofer norm) $\tilde{H_b}$ of $H_b$ for every positive real number $b$, such that $\tilde{H}_b$ is admissible of slope $b$ at infinity for all $b$, that the functions $\tilde{H}_b$ coincide on $W$ for all $b$, and that we still have $\rho(\tilde{H}_\delta,i_\delta(\alpha))>a$ (this is possible thanks to the continuity of $\rho$). We then obtain the following commutative diagram:

\[\begin{tikzcd}
    & SH^{=0}_{*,0}(W) \arrow[r,"i_\delta"]\arrow[rd,"i_a"] & SH_{*,0}^{<\delta}(W)\simeq FH_{*,0}(\tilde{H}_\delta) \arrow[r]\arrow[d] & FH_{*,0}^{\geq a}(\tilde{H}_\delta) \arrow[d,"\simeq"] \\
    & & SH_{*,0}^{<a}(W)\simeq FH_{*,0}(\tilde{H}_a) \arrow[r] & FH_{*,0}^{\geq a}(\tilde{H}_a)
\end{tikzcd}\]

Since $a>c(\alpha)$, we have $i_a(\alpha)=0$. We now use the spectral invariant to obtain a contradiction. Since $\rho(\tilde{H}_\delta,i_\delta(\alpha))>a$, the class $i_\delta(\alpha)$ is not represented entirely in $FH_*^{<a}(\tilde{H}_\delta)$, which shows that $\alpha$ is mapped to a nonzero element of $FH_*^{\geq a}(\tilde{H}_\delta)$ by the first row of the diagram. As the last vertical morphism is an isomorphism (this comes from the fact that $\tilde{H}_a$ and $\tilde{H}_\delta$ coincide on $W$, which contains all orbits of action $\geq a$ and the associated Floer trajectories), this contradicts the fact that $i_a(\alpha)=0$.

\qed

Combining Corollary \ref{c_HZ < c_s} and Proposition \ref{c_s < c_SH}, we get:

\begin{theorem}
\label{c_HZ < c_s < c_SH}
    For any relatively compact open subset $U$ of $W$, we have
    $$c_{HZ}^0(W,\overline{U};L) \leq c_\sigma(W,U;L) \leq c_{SH}(W,U;L).$$
\end{theorem}

We end by explaining how we deduce Theorem \ref{critère homologique => cordes} from Theorem \ref{c_HZ < c_s < c_SH}. Let $\overline{U}_\text{int}$ denote the complement of a collar neighborhood of $\partial U$ inside $U$. Notice that it is possible to shrink this neighborhood so that the pairs $(\overline{U}\cap L,\partial\overline{U}\cap L)$ and $(\overline{U}_\text{int}\cap L,\partial \overline{U}_\text{int}\cap L)$ are homotopy equivalent. Looking at the commutative diagram
\[\begin{tikzcd}
    & H_*(L,\partial L) \arrow[r,"j_{U!}"] & H_*(U\cap L,\partial U\cap L) \arrow[d,"\simeq"] \\
    & & H_*(U_\text{int}\cap L,\partial U_\text{int}\cap L) \arrow[lu,"j_{U_\text{int}!}",leftarrow]
\end{tikzcd}\]
we see that $\text{Ker}(j_{U!})=\text{Ker}(j_{U_{\text{int}}!})$, so 
$c_{SH}(W,U;L) = c_{SH}(W,U_{\text{int}};L).$
Applying Theorem \ref{c_HZ < c_s < c_SH} to $U_\text{int}$, we obtain that under the hypothesis of Theorem \ref{critère homologique => cordes},
$$c_{HZ}^0(W,\overline{U}_\text{int};L)\leq c_{SH}(W,U_{\text{int}};L) = c_{SH}(W,U;L)<\infty$$
which yields the almost-existence property by Theorem \ref{c(W,Z;L) => presque existence}. \qed \\

\newpage
\bibliographystyle{alpha}
\bibliography{bibliographie}

\begin{thebibliography}{BDHO24}

\bibitem[BC07]{Barraud_Cornea}
Jean-Fra{\c{c}}ois Barraud and Octav Cornea.
\newblock Lagrangian intersections and the {Serre} spectral sequence.
\newblock {\em Ann. Math. (2)}, 166(3):657--722, 2007.

\bibitem[BCS25]{ConjCordesConormal}
Filip Bro{\'c}i{\'c}, Dylan Cant, and Egor Shelukhin.
\newblock The chord conjecture for conormal bundles.
\newblock {\em Math. Ann.}, 392(3):3959--4021, 2025.

\bibitem[BDHO24]{dg-floer2}
Jean-Fran{\c{c}}ois Barraud, Mihai Damian, Vincent Humili{\`e}re, and Alexandru Oancea.
\newblock Floer {Homology} with {DG} {Coefficients}. {Applications} to cotangent bundles.
\newblock Preprint, {arXiv}:2404.07953 [math.{SG}] (2024), 2024.

\bibitem[BEE12]{BourgeoisEkholmEliashberg12}
Fr{\'e}d{\'e}ric Bourgeois, Tobias Ekholm, and Yakov Eliashberg.
\newblock Effect of {Legendrian} surgery.
\newblock {\em Geom. Topol.}, 16(1):301--389, 2012.

\bibitem[BK22]{BK22}
Gabriele Benedetti and Jungsoo Kang.
\newblock Relative {Hofer}-{Zehnder} capacity and positive symplectic homology.
\newblock In {\em Symplectic geometry. A festschrift in honour of Claude Viterbo's 60th birthday. In 2 volumes}, pages 99--130. Cham: Birkh{\"a}user, 2022.

\bibitem[CE12]{Stein2Weinstein}
Kai Cieliebak and Yakov Eliashberg.
\newblock {\em From {Stein} to {Weinstein} and back. {Symplectic} geometry of affine complex manifolds}, volume~59 of {\em Colloq. Publ., Am. Math. Soc.}
\newblock Providence, RI: American Mathematical Society (AMS), 2012.

\bibitem[FHW94]{ReliqueFloer1}
Andreas Floer, Helmut Hofer, and Krzysztof Wysocki.
\newblock Applications of symplectic homology. {I}.
\newblock {\em Math. Z.}, 217(4):577--606, 1994.

\bibitem[Flo89]{Floer_Wittenscomplex}
Andreas Floer.
\newblock Witten's complex and infinite dimensional {Morse} theory.
\newblock {\em J. Differ. Geom.}, 30(1):207--221, 1989.

\bibitem[GG04]{c_relative}
Viktor~L. Ginzburg and Ba{\c{s}}ak~Z. G{\"u}rel.
\newblock Relative {Hofer}-{Zehnder} capacity and periodic orbits in twisted cotangent bundles.
\newblock {\em Duke Math. J.}, 123(1):1--47, 2004.

\bibitem[HT13]{ConjCordesDim3}
Michael Hutchings and Clifford Taubes.
\newblock Proof of the {Arnold} chord conjecture in three dimensions. {II}.
\newblock {\em Geom. Topol.}, 17(5):2601--2688, 2013.

\bibitem[HZ94]{HoferZehnder}
Helmut Hofer and Eduard Zehnder.
\newblock {\em Symplectic invariants and {Hamiltonian} dynamics}.
\newblock Basel: Birkh{\"a}user, 1994.

\bibitem[Iri14]{Irie}
Kei Irie.
\newblock Hofer-zehnder capacity of unit disk cotangent bundles and the loop product.
\newblock {\em J. Eur. Math. Soc. (JEMS)}, 16(11):2477--2497, 2014.

\bibitem[KKK22]{KimKimKwon22}
Joontae Kim, Seongchan Kim, and Myeonggi Kwon.
\newblock Equivariant wrapped {Floer} homology and symmetric periodic {Reeb} orbits.
\newblock {\em Ergodic Theory Dyn. Syst.}, 42(5):1708--1763, 2022.

\bibitem[LR20]{c_coisotrope}
Samuel Lisi and Antonio Rieser.
\newblock Coisotropic {Hofer}-{Zehnder} capacities and non-squeezing for relative embeddings.
\newblock {\em J. Symplectic Geom.}, 18(3):819--865, 2020.

\bibitem[MS16]{McDuffSalamon17}
Dusa McDuff and Dietmar Salamon.
\newblock {\em Introduction to symplectic topology}, volume~27 of {\em Oxf. Grad. Texts Math.}
\newblock Oxford: Oxford University Press, 3rd edition edition, 2016.

\bibitem[MS18]{MurphySiegel18}
Emmy Murphy and Kyler Siegel.
\newblock Subflexible symplectic manifolds.
\newblock {\em Geom. Topol.}, 22(4):2367--2401, 2018.

\bibitem[Oan04]{Oancea_homologie_symplectique}
Alexandru Oancea.
\newblock {\em A survey of {Floer} homology for manifolds with contact type boundary or symplectic homology}, volume~7 of {\em Ensaios Mat.}
\newblock Rio de Janeiro: Sociedade Brasileira de Matem{\'a}tica, 2004.

\bibitem[Oan06]{Oan06}
Alexandru Oancea.
\newblock The {K{\"u}nneth} formula in {Floer} homology for manifolds with restricted contact type boundary.
\newblock {\em Math. Ann.}, 334(1):65--89, 2006.

\bibitem[Oh96]{Oh}
Yong-Geun Oh.
\newblock Floer cohomology, spectral sequence, and the maslov class of lagrangian embeddings.
\newblock {\em Internat. Math. Res. Notices 7}, pages 305--346, 1996.

\bibitem[Rit13]{Ritter}
Alexander~F. Ritter.
\newblock Topological quantum field theory structure on symplectic cohomology.
\newblock {\em J. Topol.}, 6(2):391--489, 2013.

\bibitem[Str90]{Struwe}
Michael Struwe.
\newblock Existence of periodic solutions of {Hamiltonian} systems on almost every energy surface.
\newblock {\em Bol. Soc. Bras. Mat., Nova S{\'e}r.}, 20(2):49--58, 1990.

\bibitem[Vit87]{Indice_Maslov_Viterbo}
Claude Viterbo.
\newblock Intersection de sous-vari{\'e}t{\'e}s lagrangiennes, fonctionnelles d'action et indice des syst{\`e}mes hamiltoniens. ({Intersection} of {Lagrangian} submanifolds, action functionals and indices of {Hamiltonian} systems).
\newblock {\em Bull. Soc. Math. Fr.}, 115:361--390, 1987.

\end{thebibliography}

\noindent
Université de Strasbourg, Institut de recherche mathématique avancée, IRMA, Strasbourg, France.

\noindent
\emph{e-mail:} antoine.rodrigues@math.unistra.fr

\end{document}